# Identifying the Preschool Home Learning Experiences That Predict Early Number Skills: Evidence From a Longitudinal Study

Elena Soto-Calvo[1], Fiona R. Simmons[1], Anne-Marie Adams[1], Hannah N. Francis[1], Hannah Patel[1], & David Giofrè[2].

[1]Liverpool John Moores University
[2]University of Genova

[8714 words, excluding abstract and references]

**Abstract**

This study examines the longitudinal relationships between home learning experiences and early number skills. The counting, number transcoding and calculation skills of 274 children were assessed in the penultimate term of preschool ($M_{age}$=4:0). Prior to these assessments, parents completed questionnaires that surveyed the frequency of the children's home learning experiences. Three types of experiences were indexed: code-focused home literacy experiences that focus on the phonological and orthographic features of language, meaning-focused home literacy experiences that focus on sharing the meaning of language and text, and home number experiences. The children's language abilities (phonological awareness and vocabulary) and nonverbal abilities (inhibitory control and nonverbal reasoning) were assessed in the final term of preschool ($M_{age}$=4:3). Their number skills were reassessed in the final term of the first year of primary school ($M_{age}$=5:3). Home letter-sound interaction experiences (interactive code-focused literacy experiences) had significant longitudinal relationships with counting and number transcoding that were independent of language and nonverbal abilities. The relationship between letter-sound interaction experiences and later counting was also independent of the autoregressive influence of baseline counting ability. We extend previous findings by demonstrating that interactive code-focused home literacy experiences in the preschool period predict growth in counting skills even when a broad range of language and cognitive abilities are controlled. Supporting parents to engage in code-focused home literacy experiences may benefit pre-schoolers' counting skills.

Early number skills are basic numerical abilities that include counting, number transcoding and calculation (Krajewski & Schneider, 2009ab). These foundational skills often emerge prior to formal schooling (e.g., LeFevre et al., 2010a; Sarnecka & Carey, 2008; Wynn, 1992), with large individual differences already apparent at school entry (Dowker, 2008; Stock et al., 2009). Children's early number skills at school entry have been repeatedly associated with later mathematics attainment (see Merkley & Ansari, 2016 for a review). Given the predictive influence of early number skills on later mathematical development, it is important to understand the language, cognitive and environmental factors that support these skills. The aim of the present study is to examine the longitudinal relationships among pre-schoolers' home learning experiences, language and cognitive abilities, and their later counting, number transcoding and calculation skills. It evaluates both within-domain relationships between home number experiences and early number skills, and cross-domain relationships between home literacy experiences and early number skills.

**Language and cognitive skills**

Within the pathways model of mathematics, LeFevre et al. (2010a) identified language abilities as one of the important precursors of early mathematical skills. Phonological skills are proposed to underpin early numeracy development by supporting the development of language-symbol associations when acquiring the names of number symbols (LeFevre et al., 2010a). It is argued that children with stronger phonological representations of numbers are better at counting and calculating because these tasks require the retrieval and usage of verbal codes (see Simmons & Singleton, 2008; Soto-Calvo et al., 2015). Phonological awareness has been associated with later mathematical attainment in young children (Barnes et al., 2011; Pupura et al., 2011; Simmons et al., 2008). Furthermore, pre-schoolers' phonological awareness predicts the development of specific early number skills, including both counting and calculation (Barnes et al., 2011; Krajewski & Schneider, 2009b; Koponen et al., 2013; LeFevre et al., 2010a; Moll et al., 2015; Soto-Calvo et al., 2015). Children's vocabulary has also been associated with their mathematical attainment (Barnes et al., 2011; Romano et al., 2010). Vocabulary is proposed to support mathematical development through the acquisition and understanding of mathematical vocabulary (Moll et al., 2015; Purpura & Logan, 2015; Purpura & Reid, 2016; Toll & Van Luit, 2014). Pre-schoolers' vocabulary skills are predictive of their early number skills including counting and number transcoding (Moll et al., 2015; Purpura & Ganley, 2014). Combined this research suggests that phonological and vocabulary skills support the development of verbally mediated aspects of early numeracy.

Executive functioning (EF) has also been identified as a precursor of early numeracy (see Raghubar & Barnes, 2017 for a review). EF is an overarching term for the cognitive processes that control goal-directed behaviour. These processes include the updating of working memory, inhibitory control (IC) and the shifting of the attentional focus (Wiebe et al., 2011). There are a number of ways that EF may support early number skills. For example, counting objects requires updating of working memory to maintain an accurate representation of the changing count sequence and inhibiting previously counted objects so they are not counted twice. Similarly, updating of working memory is required during calculation as the quantities transform. Pre-schoolers' performance on EF tasks relates to their concurrent and longitudinal performance on standardised mathematics tests (e.g., Bull et al., 2011; McClelland et al., 2013), with several studies suggesting these relationships are causal (e.g., Blair et al., 2015; Cameron et al., 2012; Davidse et al., 2015; Matthews et al., 2009; McClelland et al., 2014). Researchers have also reported positive links between EF measures and specific early number skills including counting, number transcoding and calculation at preschool and kindergarten ages (Hornung et al., 2014; Kroesbergen et al., 2014; Lan et al., 2011; Moll et al., 2015; Purpura & Ganley, 2014; Purpura et al., 2017).

Reasoning abilities have been shown to be an important predictor of mathematics attainment in secondary and older primary-aged children (e.g., Allen et al., 2019; Caviola et al., 2014; Deary et al., 2007; Donolato et al., 2019), although their relationship with specific early number skills in young children is more nuanced. Relationships between nonverbal reasoning and measures of early number skills (that comprised mainly counting tasks) have been identified in kindergarteners (Hornung et al., 2014; Kroesbergen et al., 2014). However, these relationships were not independent of children's EF. Furthermore, although nonverbal reasoning was not independently related to a counting-focused kindergarten assessment, it did have an independent relationship with number line estimation and arithmetic in 1st grade (Hornung et al., 2014). This suggests that nonverbal reasoning may have a more limited role in early counting but a stronger role in later mathematical skills.

## The Home Learning Environment

The Home Learning Environment (HLE) is a broad term that encompasses the frequency of home learning experiences, the availability of resources that promote learning and parental attitudes towards learning. The HLE in early years is a longitudinal predictor of academic outcomes including mathematics (Anders et al., 2013; LeFevre et al., 2009; Melhuish et al., 2008a; Skwarchuk et al., 2014). Research adopting broad indices of the HLE suggest that a higher quality HLE is associated with stronger numeracy (e.g., Anders et al., 2012; Melhuish et al., 2008b). Examining which aspects of the HLE exert an influence on early number skills development is particularly important because intervention studies have indicated that the preschool HLE can be modified (Niklas et al., 2016; Suskind et al., 2015).

### *The Home Literacy Environment*

Home literacy models provide an analysis of the different types of home literacy experiences that children experience (Phillips & Lonigan, 2009; Sénéchal & LeFevre, 2002). Within these models experiences are classified into those that focus on phonological and orthographic aspects of language (termed 'formal', 'inside-out' or 'code-focused'), and those that prioritise enjoyment and focus on sharing the meaning of oral or written language (termed 'informal', 'outside-in' or 'meaning-focused'). Code-focused home literacy experiences would include teaching children about letters and words, whereas the quintessential meaning-focused home literacy experience is shared reading. The home literacy model (Sénéchal & LeFevre, 2002) proposes that code-focused experiences support the development of emergent literacy skills (including alphabetic knowledge). In contrast, shared reading is proposed to support semantic oral language skills. Empirical studies have largely supported the predictions of this model. Meaning-focused experiences are associated with vocabulary and wider semantic skills (Frijters et al., 2000; Hamilton et al., 2016; Hood et al., 2008; Manolitsis et al., 2013; Sénéchal, 2006; Sénéchal & LeFevre, 2002; Sénéchal et al., 2008; Skwarchuk et al., 2014, but cf. Evans et al., 2000). Code-focused experiences are associated with alphabetic knowledge and word reading (e.g., Hood et al., 2008; Sénéchal & LeFevre, 2002, 2014). In some studies code-focused experiences have also been related to phonological skills (Foy & Mann, 2003; Hamilton et al., 2016), but this finding has not been consistently replicated (Hood et al., 2008; Napoli & Purpura, 2018).

Positive links have been identified between home literacy experiences and early numeracy (Anders et al., 2012; Baker, 2014; Napoli & Purpura, 2018; Soto-Calvo et al., 2020), although this relationship is not consistently found (Huntsinger et al., 2016; LeFevre et al., 2009; Segers et al., 2015). Manolitsis et al. (2013) examined the relationships between code- and meaning-focused home literacy experiences and counting knowledge using path analyses. Code-focused experiences were related to counting knowledge whereas meaning-focused experiences were not. Although not examined within their regression analyses, the correlations reported by Napoli and Purpura (2018) show a similar pattern with code-focused experiences having stronger relationships with early mathematics than meaning-focused experiences. It is possible that previous studies have failed to identify relationships between home literacy indices and numeracy outcomes because the home literacy indices employed focused on meaning- rather than code-focused items.

### ***The Home Numeracy Environment***

A model of the home numeracy environment that echoes the formal and informal distinction within the home literacy model has been developed (Skwarchuk et al., 2014). Within this model, home number experiences that involve numerical skills explicitly are considered ‘formal’ or ‘direct’. Experiences that focus on enjoyment or play, where the number-related content is more incidental, are considered ‘informal’ or ‘indirect’. Typical ‘formal’ experiences would include encouraging a child to complete calculations or teaching them to count objects. ‘Informal’ experiences would include playing shops or playing board games with a dice. Alongside the formal/informal distinction many studies also classify number experiences as basic or advanced, the latter being more stretching for the age-group of children studied (Elliott & Bachman, 2017; Thompson et al., 2017). Additionally, some studies consider mathematical experiences beyond numeracy e.g., spatial and patterning experiences (Hart et al., 2016; Huntsinger et al., 2016; Zippert & Rittle-Johnson, 2020). Many studies have shown that *some* measures of the home numeracy environment are related to numerical attainment (Anders et al., 2012; del Río et al., 2017; Hart et al., 2016; Huntsinger et al., 2016; Napoli & Purpura, 2018; Niklas & Schneider, 2013; Skwarchuk et al., 2014; Sonnenschein et al., 2016), although others have failed to identify any significant relationships (Blevins-Knabe et al., 2000; Lehrl et al., 2019; Missall et al., 2015). These differences are most likely underpinned by how the home numeracy environment is conceptualised and assessed. Patterns consistent with this interpretation can be identified. First, home number experiences categorised as ‘formal’ or ‘direct’ appear stronger predictors of symbolic number skills than those classified as ‘informal’ or ‘indirect’ (Huntsinger et al., 2016; Skwarchuk et al., 2014; Thompson et al., 2017). Second, advanced home number experiences bear stronger relationships with numeracy outcomes than basic experiences (del Río et al., 2017; LeFevre et al., 2010b; Ramani et al., 2015; Skwarchuk, 2009; Skwarchuk et al., 2014; Thompson et al., 2017). Overall, this research suggests that more frequent home experiences that are challenging and focus explicitly on number *may* support the development of children's number skills by increasing their exposure to, and practice of, number-related content.

**The Present Study**

Within the present study, the relationships between home number experiences, code- and meaning-focused home literacy experiences and early number skills were explored. The extent that any cross-domain relationships between home literacy experiences and children's early number skills are direct or indirect via language and nonverbal abilities was analysed. Previous studies that have analysed the relative impact of code- and meaning-focused literacy experiences on early mathematics and number skills (Manolitsis et al., 2013; Napoli & Purpura, 2018) have not considered the extent that such relationships are underpinned by children's language abilities. Furthermore, we examine the within- and cross-domain influences of home learning experiences on three specific number skills, counting, number transcoding and calculation, using dual measures of each number skill. This extends previous research by examining the *differential* impact of different aspects of the HLE on distinct early number skills.

We focused on counting, number transcoding and calculation skills because they have been consistently associated with children's later mathematical attainment (e.g., Krajewski & Schneider, 2009a,b; Moll et al., 2015; Passolunghi et al., 2007; Purpura et al., 2011; Soto-Calvo et al., 2015; Stock et al., 2009; van Marle et al., 2014). They are also core components of the Early Years Curriculum in England (Department for Education, 2013, 2017; Testing and Standards Agency, 2017). Theoretical models of early numeracy propose that counting, number transcoding and calculation skills are distinct early number skills that may be supported by different cognitive abilities (Krajewski & Schneider, 2009b; LeFevre et al., 2010a).

It is particularly important to understand the factors influencing children's ability to use counting to enumerate sets effectively as such cardinal counting skills are considered an essential building block in numeracy development (Sarnecka & Carey, 2008; Wynn, 1992). Krajewski and Schneider's (2009b) model of early number skills views cardinal counting as a pre-requisite for the development of higher order skills such as symbolic number magnitude comparisons and calculations. The importance of cardinal counting skills is also supported by evidence that the establishment of cardinal counting principles and procedures is associated both with an acceleration of the understanding of the magnitude of numbers (Geary et al., 2018), and an understanding of the principles underlying addition (Sarnecka & Carey, 2008).

In summary, we extend previous research by examining the relationships among number, code-focused literacy and meaning-focused literacy experiences and three core early number skills using a methodology that enabled the influence of a broad range of language and cognitive skills to be evaluated. We addressed the following research questions:

1. *Does the frequency of home number and home literacy experiences independently predict later counting, number transcoding and calculation skills?* In view of research reporting positive links between formal or direct home number experiences and symbolic aspects of mathematics (Anders et al., 2012; del Río et al., 2017; Hart et al., 2016; Huntsinger et al., 2016; Napoli & Purpura, 2018; Niklas & Schneider, 2013; Skwarchuk et al., 2014; Sonnenschein et al., 2016), we hypothesised that the frequency of home number experiences would explain unique variance in children's counting, number transcoding and calculation skills. Given previous findings of stronger relationships between code-, rather than meaning-focused, home literacy experiences and numeracy skills in young children (Manolitsis et al., 2013; Napoli & Purpura, 2018), we hypothesised that the frequency of code-, but not meaning-focused literacy experiences, would explain unique variance in later counting, number transcoding and calculation skills.

2. *To what extent do language and nonverbal abilities predict later counting, number transcoding and calculation skills? Are any relationships between the frequency of home literacy experiences and early number skills direct or indirect via children's language and nonverbal abilities?* Within the current study, we administered a broad range of language measures that assessed both phonological awareness and vocabulary, and nonverbal abilities that assessed both non-verbal reasoning and inhibitory control (IC). Given the need to restrict the assessment battery administered to pre-schoolers, we could not assess all aspects of EF. IC skills were chosen since these develop earlier than other aspects of EF (Best & Miller, 2010; Garon et al., 2008; Hughes, 2011) and in young children have stronger relationships with maths performance than other aspects of EF (Blair et al., 2015; Cameron et al., 2012; Davidse et al., 2015; Jacob & Parkinson, 2015; Matthews et al., 2009; McClelland et al., 2014).

Consistent with the results of previous studies (Barnes et al., 2011; Krajewski & Schneider, 2009b; Koponen et al., 2013; Lan et al., 2011; LeFevre et al., 2010a; Moll et al., 2015; Purpura & Ganley, 2014; Purpura et al., 2017; Soto-Calvo et al., 2015), we hypothesised that both language and nonverbal abilities would explain significant variance in counting, number transcoding and calculation skills. Furthermore, given previously identified relationships between language abilities and both home literacy experiences (Frijters et al., 2000; Hamilton et al., 2016; Hood et al., 2008; Manolitsis et al., 2013; Sénéchal, 2006; Sénéchal & LeFevre, 2002; Skwarchuk et al., 2014) and early number skills (Moll et al., 2015; Soto-Calvo et al., 2015), we anticipated that the relationships between home literacy experiences and early counting, number transcoding and calculation would be at least partially indirect via children's language abilities.

2. *Can home literacy and/or number experiences predict growth in counting, number transcoding and calculation skills?* We determined if home experiences could predict *growth* in early number skills by examining whether longitudinal relationships between home learning experiences and later number skills persisted once the autoregressive effect of children's baseline early number skills were controlled.

## Method

### Participants

Parent-child dyads were recruited in children's final year of preschool (the academic year in which they turned four). There were 274 parents (254 females) and 274 children (146 females). The children were assessed at three time points; spring term of preschool year ($T_1$, $M_{age}$ = 4:0 *SD*=3.63 months), summer term of preschool year ($T_2$, $M_{age}$ = 4:3 *SD*=3.67 months) and after they had transferred to primary school in the summer term of Reception Year ($T_3$, $M_{age}$ =5:3, *SD*=3.64 months). All early years providers in England (which include all preschool settings and Reception classes) must follow the Early Years Foundation Stage Framework (Department for Education, 2017). This framework stipulates both broad areas of learning and development that practitioners must focus on, as well as safeguarding and welfare requirements.

The parents completed a questionnaire indexing demographic factors and home learning experiences prior to their children completing the assessments at $T_1$. From the original 274 children recruited at $T_1$, 265 (96.72%) were retained at $T_2$ and 232 (84.67%) were retained at $T_3$. Attrition was higher from $T_2$ to $T_3$ because the children moved from their preschool settings to primary schools. The primary schools into which 26 (4.01%) children had transferred did not consent for the study to continue and the primary school placements of 10 (3.65%) children could not be traced. Additionally, four children (1.46%) moved away, two (0.73%) were persistently absent and one (0.36%) withdrew assent.

The postcode deprivation decile for each household was obtained from the English indices of deprivation 2015 online open data of the United Kingdom (Department for Communities and Local Government, http://imd-by-postcode.opendatacommunities.org/). The mean deprivation level was close to the national average (*M*=5.42, *SD*=3.32). Three respondents did not supply their postcode. Parental qualifications were coded according to the UK National Qualification Framework (https://www.gov.uk/what-different-qualification-levels-mean/list-of-qualification-levels). This scale levels qualifications from 1 (qualifications equivalent to a lower grade GCSE, typically taken by 16-year-olds) to 8 (doctoral level qualifications). Parental highest level of education was diverse, with a mean which was broadly equivalent to two years of post-secondary education (*M*=4.75, *SD*=2.00). Four respondents did not report their qualifications.

**Measures**

***Questionnaire***

**Home experiences.** Parents reported the frequency on a 6-point Likert scale ranging from *never* to *several times a day* that their child experienced 32 activities at home. There were eight number experiences, eight meaning-focused literacy experiences and seven code-focused literacy experiences. In addition, there were nine filler items that were not analysed (e.g., rides a scooter, balance bike or bike). The different types of items were randomly ordered within the questionnaire. Three factor analyses were conducted to assess whether the items relating to number, meaning-focused literacy and code-focused literacy experiences formed reliable and internally consistent scales. Principal Axis Factoring with a Promax rotation and Kaiser normalisation was utilised including only items deemed statistically reliable (i.e., did not suffer from lack of variability in response and were largely correlated with the remaining items in the scale). Missing item responses were replaced with the item mean. Appendices A to C show the key descriptive statistics for the included and excluded items, and the inter-item reliability of the scales.

**Book exposure.** Parents were presented with 21 potential preschool book titles and asked to indicate which titles they believed were real. They were given the response options 'real', 'made up' and 'don't know'. The instructions encouraged them not to guess. Of the 21 book titles, 15 were real and six were made up. All of the real titles were fiction books suitable for preschoolers. Parental responses to the real and fictional titles are shown in Appendix D. This book exposure index is based on similar book exposure indices that have been used successfully to assess shared reading in previous studies (e.g., Hamilton et al., 2016; Puglisi et al., 2017; Skwarchuk et al., 2014). The same formula as Skwarchuk et al. (2014) was utilised [(Storybook titles correctly identified - Foils identified as real books) / total number of actual books x 100] to reduce any residual influence of guessing.

***Child Assessments***

**Language Assessments.** Two subtests from the Preschool and Primary Inventory of Phonological Awareness (PIPA, Dodd et al., 2000) were administered to assess phonological awareness. In *Alliteration Awareness,* the child had to identify the word from a set of four that started with a different sound. In *Rhyme Awareness,* the child had to identify the word from a set of four that did not rhyme with the others. Both tests comprise 2 practice items and 12 experimental items. Vocabulary was assessed with two standardised measures. In the *Naming Vocabulary* subtest from the British Ability Scales III (BAS-3, Elliott & Smith, 2011) the child had to name pictures presented to them. In the *Receptive Vocabulary* subtest from the Wechsler Preschool and Primary Scale of Intelligence-Fourth UK Edition (WIPPSI-IV-UK, Wechsler, 2013) the child had to point at the picture from a set of four that best matched the word said by the researcher.

**Cognitive Assessments.** IC was assessed with two experimental computerised tasks previously used with pre-schoolers. In the *Fish/Shark* task (Wiebe et al., 2011) the child had to press a key when shown a fish (75% of the trials) but inhibit this response when shown a shark (25% of the trials). The *d'* index was calculated by subtracting the *z*-score value of the hit rate right-tail *p* value from the *z*-score value of the false alarm rate right-tail *p* value (Macmillan & Creelman, 2005). This sensitivity index represents how accurately the child detects fishes and rejects sharks. In the *Big/Little Stroop* task (modified from Kochanska et al., 2000) the child saw a large outline of an animal with smaller animal outlines presented within it. The large outline appeared briefly first for 750ms. The child had to name the smaller animals within the outline. The trials were equally split between congruent trials (where the outline animal matched the smaller ones within it) and incongruent trials (where the outline animal differed from the small ones within it). The number of correct responses to the incongruent trials was recorded. We chose our two IC measures because they both capture variance in this age group effectively (Clark et al., 2014; Kochanska et al., 2000; Wiebe et al. 2011). Nonverbal reasoning was assessed with two standardised measures from the British Ability Scales III (BAS-3, Elliott & Smith, 2011). In the *Matrices* subtest, the child had to select the picture from a set of four that best completed a four-picture pattern. In the *Picture Similarities* subtest, the child had to place a card under a picture from a choice of four that best matched the picture on the card.

**Early Number Skills Assessments.** The counting, number transcoding and calculation assessments were developed for this longitudinal study, further details including the individual items are given in Soto-Calvo et al. (2020). With the exception of sequential counting, all the number assessments began with one practice item for which feedback was provided. One point was awarded for each correct response in the subsequent experimental items. To reduce administration time and stress (due to the presentation of numerous items beyond a child's capabilities) stop rules were applied at both time points and start rules were applied at $T_3$. At $T_1$ and $T_3$, the child was stopped after three consecutive incorrect answers or after three or more incorrect answers within a block. If a child completed more than one block without error at $T_1$, at $T_3$ the child started the task at the beginning of the last block in which all items were answered correctly at $T_1$. Credit was given for preceding items not administered. If starting on a later block at $T_3$, the previous block was administered if the child failed any items on the first block administered.

***Sequential Counting.*** Children were asked to count aloud to a cuddly toy starting from one to as high as they could. The highest number recited in the correct order was recorded.

***Give Me X.*** The child was given a set of objects (e.g., 10 pigs) and asked to select a subset (e.g., 5 pigs). In each item, the child was asked to place a specific number of toy animals on a drawing of a farm or a house (e.g., "Can you put two ducks in the pond?"). The task consisted of three blocks of five items. Items in the first block were magnitudes below 10 and the child had to select the items from a box containing 10 items, items in the second block were numbers from 10 to 20 and items in the third block were numbers from 20 to 30. In these latter two blocks, the child had to select the items from a box containing 35 items.

***Counting Objects.*** The child was asked how many animal pictures were presented on a card (e.g., "How many birds are there?"). There were 20 cards with pseudo-randomly distributed pictures of the same animal on each card. The last number-word spoken was recorded as their answer. The cards were grouped into four blocks each consisting of five items. The first block presented quantities below 10, the second presented quantities from 10 to 19, the third presented quantities from 20 to 30 and the fourth presented quantities ranging from 35 to 97.

***Numeral Recognition.*** The researcher presented a card displaying nine different numerals in a pseudo-random arrangement and asked the child to point to a specific number (e.g., “Can you point to number five?”). This task consisted of four blocks at $T_1$ and five blocks at $T_3$ to prevent ceiling effects. Each block consisted of five items. In the first block, the card presented the nine single-digit numerals. In the second and third blocks, the card presented two-digit numerals selected from the range 11 to 19 and from 20 to 90, respectively. In the fourth and fifth blocks, the card presented three-digit numerals selected from the range 100 and 200 and from 200 to 999, respectively.

***Numeral Reading.*** The researcher pointed to a printed numeral on a card containing five numerals and asked the child to name it. This task consisted of four blocks each containing five items. The first block presented numerals selected from the range 1 to 9. The second and third blocks presented two-digit numerals from the range 10 to 19 and 20 to 99, respectively. The fourth block presented three-digit numerals selected from the range 100 to 199.

***Additions and Subtractions.*** The researcher presented addition and subtraction problems to the child in the form of a story (e.g., “If you put two horses on the path and you add one more, how many horses would there be?”). Animal toys and a drawing of a farm or a house were available for the child to use to support their calculation. The child was asked to provide a verbal response with the last number-word spoken recorded as their answer. Both the addition and subtraction sections consisted of three blocks each containing four items. In the first block, the child had to add or subtract one or two to an addend or minuend below five. In the second block, the child had to add or subtract two or three to an addend or minuend of four, five or six. In the third block, for addition, the child had to add three or four to a single-digit addend with the answer always being above 10. For subtraction, the child had to subtract three or four from a minuend above 10.

**Mathematics Assessment.** The Early Number Concepts test from the British Ability Scales III (BAS-3 ENC, Elliot & Smith, 2011) was administered. This test assesses a broad range of age-appropriate mathematical concepts. Internal consistency for this test is .84.

**Procedure**

The data presented here is from a wider longitudinal study which examines the relationships among environmental and cognitive factors, and academic attainment (see Soto-Calvo et al., 2020 for concurrent analyses of the home learning experiences and early number skills and further details of the sample and tasks). Ethical approval was granted by the Liverpool John Moores University research ethics panel in August 2016. Written consent was gained from the educational settings' managers. The children completed individual assessment sessions in a quiet area of their preschool ($T_1$ and $T_2$) and school ($T_3$). Each session lasted approximately 15 minutes. At $T_1$ and $T_{3,}$ the children completed the early number skills tasks in two sessions. The first session consisted of numeral recognition, give me X and additions. The second session consisted of numeral reading, counting objects, subtractions and sequential counting. At $T_{3,}$ the children also completed the standardised mathematics measure in an additional third session. At $T_{2,}$ the children completed the language and cognitive assessments in three sessions. In the first session, the children completed the Matrices, Naming Vocabulary and Alliteration Awareness subtests. In the second session, they completed the Receptive Vocabulary and Picture Similarities subtests. In the third session, they completed the Big/Little Stroop task, the Rhyme Awareness subtest and the Fish/Shark task. For all standardised tests the standardised administration, progression and scoring rules were applied.

**Analysis**

We conducted preliminary analyses before addressing the research questions. Descriptive statistics and correlations between the individual measures were obtained. The structure of the early number skills, language and cognitive variables was assessed using Confirmatory Factor Analysis (CFA). This was to ensure that the measures were grouped appropriately when we constructed composite variables for the longitudinal analyses. At this stage, the validity of the early number skills measures was also assessed by examining the concurrent relationships between the early number skills and the standardised mathematics measure.

Once the structure of the early number skills, language and cognitive variables was confirmed, we created composite variables that indexed the core number, language and cognitive factors. To address our research questions three path analyses were conducted. Model 1 addressed the first research question and assessed whether code-focused home literacy, meaning-focused home literacy, and home number experiences at $T_1$ predicted children's early number skills at $T_3$. Model 2 addressed the second research question by assessing whether the aspects of the HLE that predicted early number skills at $T_3$ in Model 1 continued to predict children's early number skills when the language and cognitive variables at $T_2$ were included in the model. In this model, direct longitudinal relationships between the home learning experiences at $T_1$ and the early number skills at $T_3$ were assessed alongside indirect relationships via the language and cognitive variables at $T_2$. Finally, Model 3 addressed the third research question by expanding Model 2 to include children's initial number skills at $T_1$ as predictors of the number skills at $T_3$. This inclusion of initial skill level in the model controls for autoregressive effects and enables the extent that home learning experiences can predict growth in early number skills to be assessed.

The R program (R Core Team, 2014) with the "lavaan" library (Rosseel, 2012) was used to conduct the CFAs and the path analyses. Model fit was assessed using various indices according to the criteria suggested by Hu and Bentler (1999). We considered the chi-square ($\chi^2$), the comparative fit index (*CFI*), the non-normed fit index (*NNFI*), the standardised root mean square residual (*SRMR*) and the root mean square error of approximation (*RMSEA*) to evaluate model fit. Only participants with data for all variables analysed were included in the CFAs and path analyses. The *N* for each CFA and path analysis is shown alongside the figures. Participants were excluded from the path analyses if they had missing data. We also repeated the analyses using the Full Information Maximum Likelihood method for estimating missing data. The two approaches resulted in minimal differences and therefore we maintained the initial analyses (see figure notes for details of any differences when FIML was applied).

## Results

### Preliminary Analyses: The Early Number Skills

Table 1 shows correlations between the demographic variables, the home learning scales and the early number skills at $T_1$ and at $T_3$ together with descriptive statistics for these variables. We performed a series of CFAs to confirm the factorial structure of the early number skills at $T_1$ and at $T_3$. In our previous analysis of the preschool number skills (Soto-Calvo et al., 2020) a three-factor structure provided a good fit for the data. The three factors were counting (counting objects and give me X), number transcoding (numeral recognition and numeral reading) and calculation (additions and subtractions). The sequential counting task was excluded because it reduced the fit of the model whether included as a separate factor or included with the other counting measures (Soto-Calvo et al., 2020). At $T_1$ the fit of this model was good, $\chi^2(6)=16.40$, $p=.012$, $RMSEA=.09$, $SRMR=.03$, $CFI=.99$, $NNFI=.96$. We tested the same model with the $T_3$ early number skills. The fit of the model was completely satisfactory, $\chi^2(6)=2.87$, $p=.83$, $RMSEA=.00$, $SRMR=.01$, $CFI=1.00$, $NNFI=1.01$.

Performance on the BAS-3 ENC is comparable to that of the standardisation sample (present sample $M=49.80$, $SD=11.67$, standardisation sample $M=50$, $SD=10$), suggesting that the mathematical attainment of our sample is broadly consistent with that of children in the UK. A single-step multiple linear regression was conducted to test the early number skills' criterion validity. Counting, number transcoding and calculation skills at $T_3$ were entered simultaneously and explained a significant 43% of the variance in the BAS-3 ENC ($r^2=.43$, $F(3, 227)=58.12$, $p < .001$; $b_{counting\ T3}=.20$, $p < .01$, $b_{number\ transcoding\ T3}=.26$, $p < .001$, and $b_{calculation\ T3}=.30$, $p < .001$,

*N*=231). Our counting, number transcoding and calculation factors, although related, are separable and together explain a significant proportion of variance in mathematical attainment.

**Preliminary Analyses: The Language and Cognitive Skills**

Table 2 shows the means, standard deviations and correlations between the cognitive and language measures obtained at $T_2$. We ran two CFAs to examine the factorial structure of the language and cognitive variables. In the first CFA, we organised the measures into four factors, phonological awareness (Rhyme awareness and Alliteration awareness), vocabulary (Naming vocabulary and Receptive vocabulary), IC (Fish/Shark d' and Big/Little Stroop) and nonverbal reasoning (Matrices and Picture similarities). The fit of this model was not strong, $\chi^2(14)=16.03$, $p=.31$, $RMSEA=.03$, $SRMR=.03$, $CFI=.99$, $NNFI=.99$ so we tested a further model where the measures were collapsed into two factors. In this second CFA, we created a language factor consisting of the phonological and vocabulary measures, and a nonverbal factor consisting of the IC and the nonverbal reasoning measures. This two-factor structure had a better fit, $\chi^2(19)=31.89$, $p=.03$, $RMSEA=.05$, $SRMR=.05$, $CFI=.96$, $NNFI=.94$ and was therefore retained for further scrutiny.

**Path Analyses**

We used path analyses to examine the longitudinal relationships between home learning experiences, language and nonverbal abilities, and the early number skills following a methodology similar to Manolitsis et al. (2013) and Moll et al. (2015). In order to conduct the path analyses we created composite scores consistent with the structures confirmed by the CFAs. These composites were created from the mean of the z-scores of the component variables. At both $T_1$ and $T_3$ counting was created from give me x and counting objects, number transcoding from numeral reading and numeral recognition, and calculation from addition and subtraction. The language composite was created from Alliteration awareness, Rhyme awareness, Receptive vocabulary and Naming vocabulary. The nonverbal composite was created from Matrices, Picture similarities, Big/Little Stroop and Fish/Shark d'. Correlations between the early number, language and cognitive composites, the demographic variables and the home learning scales are shown in Table 3.

***Assessing the Longitudinal Relationships Between the HLE and Early Number Skills***

Model 1 addressed research question 1 and examined the extent that home number and home literacy experiences during preschool predicted later early number skills. It predicted counting, number transcoding and calculation skills at $T_3$ from the background variables (postcode deprivation decile and parental qualifications), the home number experiences scale and the letter-sound interaction experiences scale. We did not include the meaning-focused literacy experiences or book exposure scales as potential predictors because the raw correlations between these indices and the number skills composites were non-significant and close to zero (see Table 3). To avoid collinearity issues we chose to include the letter-sound interaction experiences scale as the code-focused home literacy index because it had stronger relationships with the early number composites than letter activities (see Table 3). Model 1 (illustrated in Figure 1) was saturated. The fit was completely adequate. Only letter-sound interaction experiences significantly predicted the three early number skills at $T_3$. The proportion of variance explained in the early number skills was modest (Counting 11%, Number transcoding 12% and Calculation 8%).

***Assessing the Inter-relationships Between Code-focused Home Literacy Experiences, Language and Nonverbal Abilities, and Early Number Skills***

Model 2 addressed research question 2. It examined the extent that language and nonverbal abilities predicted the early number skills to determine whether the relationships between letter-sound interaction experiences and early number skills (identified in Model 1) were direct or indirect via children's language and cognitive abilities. We tested both direct paths, from letter-sound interaction experiences to each early number skill and indirect paths, via both language and nonverbal abilities. This model was saturated. All the paths were statistically significant, except for letter-sound interaction experiences, which was not a statistically significant predictor of calculation at $T_3$. This path was therefore dropped and a revised model tested. The fit of the final version of Model 2 was perfectly appropriated, $\chi^2(1)=1.05$, $p=.31$, *RMSEA*=.02, *SRMR*=.02, *CFI*=1.00, *NNFI*=1.00, and all paths were significant (illustrated in Figure 2). Within this model letter-sound interaction experiences were related to both nonverbal abilities and language which in turn were related to all three number skills. In addition to these indirect relationships, letter-sound interaction experiences also directly explained significant variance in counting and number transcoding. The proportion of the variance explained in early number skills was large (Counting 25%, Number transcoding 25% and Calculation 20%).

***Assessing if Code-focused Home Literacy Experiences Predict Growth in Early Number Skills***

Finally, we used a third model to address research question 3 and explore the autoregressive effects of the children's preschool early number skills on the identified relationships. In Model 3, the relationships within Model 2 were tested with the effects of number skills at $T_1$ controlled. It should be noted that the fit of Model 3 was poor, $\chi^2(13)=78.40$, $p<.001$, *RMSEA*=.16, *SRMR*=.16, *CFI*=.87, *NNFI*=0.71, and therefore one should be cautious in its interpretation (illustrated in Figure 3). Despite the poor fit, all previous significant paths remained significant except for nonverbal abilities that no longer predicted counting or calculation at $T_3$, and letter-sound interaction experiences that just missed the traditional level of significance ($p=.052$) when predicting number transcoding at $T_3$. The proportion of variance explained in early number skills was large (Counting 26%, Number transcoding 37% and Calculation 24%).

Given that overall the model had a poor fit, we also used regressions to explore the impact of the autoregressive effects of children's baseline number skills on later number skills and address research question 3. Hierarchical linear regressions (see Table 4) were conducted with children's early number skills at $T_1$ (the autoregressor) entered at step 1, the language and nonverbal abilities entered at step 2, and the frequency of letter-sound interaction experiences entered at step 3. The results were consistent with Model 3. Letter-sound interaction experiences explained unique variance in counting skills even when the autoregressive effect of initial abilities, and language and nonverbal abilities were accounted for, but could not explain variance in number transcoding and calculation in the same model. Consequently, we can conclude that letter-sound interaction experiences predict growth in counting.

## Discussion

This study examined the extent that the frequency of home number and home literacy experiences during the preschool period predict later number skills. It extends previous research that has examined the influence of the home literacy environment on early mathematics (Anders et al., 2012; Baker, 2014; Manolitsis et al., 2013; Napoli & Purpura, 2018) by exploring the relative influence of code- and meaning-focused home literacy experiences on three specific number skills: counting, number transcoding and calculation. When a broad range of language and cognitive abilities were accounted for, letter-sound interaction experiences (an index of code-focused home literacy) directly predicted later counting and number transcoding, although not calculation. Once the autoregressive effects of initial number skill levels were controlled, letter-sound interaction experiences continued to predict *growth* in counting. In contrast, neither the frequency of number experiences nor meaning-focused home literacy experiences were independent predictors of any of the early number skills studied.

### Language, Nonverbal Abilities, and Early Number Skills

The significant pathways from language abilities to early number skills in Models 2 and 3 suggest that children’s preschool language abilities support the development of counting, number transcoding and calculation. We therefore accept the hypothesis that both preschool language and nonverbal abilities explain unique variance in counting, number transcoding and calculation in the first year of primary school. This concurs with the prediction of the pathways model (LeFevre et al., 2010a) that language is an important predictor of early number skills. This finding is also in line with previous studies reporting links between language abilities and children’s counting, number transcoding and calculation skills (Barnes et al., 2011; Koponen et al., 2013; Krajewski & Schneider, 2009b; LeFevre et al., 2010a; Moll et al., 2015; Soto-Calvo et al., 2015). Within the current study, nonverbal abilities predicted all three early number skills in Model 2 (although the variance explained was smaller than that attributed to language abilities). Of the paths from nonverbal ability to the early number skills studied, however, only the path with number transcoding remained significant in Model 3 when the autoregressive effects of preschool number skills were accounted for. These findings suggest that language abilities have a stronger impact on the development of the early number skills examined than nonverbal abilities. This pattern of findings is consistent with Moll et al. (2015) who reported that language abilities predicted a greater proportion of variance in children's early counting and number transcoding than EF. These findings are likely to be due to the verbal nature of the aspects of early number skills studied. Nonverbal reasoning and EF may have a stronger impact on other aspects of mathematics such as understanding the magnitude of symbolic numbers (see Hornung et al., 2014).

**Home Number Experiences and Early Number Skills**

Although the frequency of code-focused literacy experiences (letter-sound interaction experiences) was an independent predictor of early number skills, the frequency of home number experiences did not *independently* predict later early number skills. Therefore, we cannot accept the hypothesis that home number experiences explain independent variance in early number skills. Our home number scale comprised experiences with an explicit numerical component because these have been more consistently associated with symbolic number skills (LeFevre et al., 2010b; Swarchuk et al., 2014) and mathematical attainment (Huntsinger et al., 2016) than informal or indirect number experiences. Home number experiences were consistently related to counting, number transcoding and calculation at both $T_1$ and $T_3$. However, Model 1 demonstrated that these relationships were not independent of home letter-sound interaction experiences. The results of the current study differ from previous findings reporting positive relationships between the home numeracy environment and young children's numeracy skills (Anders et al., 2012; del Río et al., 2017; Hart et al., 2016; Huntsinger et al., 2016; Skwarchuk et al., 2014; Sonnenschein et al., 2016). Our results may differ due to the level of challenge of the number experiences we surveyed. More advanced home number experiences have a stronger relationship with symbolic number skills than more basic home number experiences (del Río et al., 2017; LeFevre et al., 2010b; Ramani et al., 2015; Skwarchuk, 2009; Skwarchuk et al., 2014; Thompson et al., 2017). A stronger relationship between home number experiences and early number skills may therefore be identified if more challenging home number experiences than the ones used in the present study are surveyed (see Elliott & Bachman, 2017 and Thompson et al., 2017 for discussions of this issue).

**Home Literacy Experiences and Early Number Skills**

We found that the frequency of code-focused home literacy experiences (indexed by the frequency of letter-sound interaction experiences) predicted children's counting, number transcoding and calculation skills. Our findings are consistent with previous evidence of positive relationships between home literacy experiences and early numeracy (Anders et al., 2012; LeFevre et al., 2010b; LeFevre et al., 2009; Melhuish et al., 2008b; Skwarchuk et al., 2014). Furthermore, our findings clarify the type of home literacy experiences that predict number skills at this early stage of development. In contrast to the relationships identified between letter-sound interaction experiences and early number skills, the frequency of meaning-focused experiences and book exposure (that indexes shared reading) did not explain significant variance in early number skills. This pattern of results where experiences that focus on orthography and phonology and the links between them, but not experiences that focus on sharing the meaning of oral or written language, are related to children's early number skills is consistent with Manolitsis et al. (2013) and Napoli and Purpura (2018). We also extend these findings by demonstrating that code-focused home literacy experiences are more robust predictors of early counting than calculation. Our evidence of a stronger association with code- rather than meaning-focused home literacy experiences and number skills could help to explain previous null findings where no association between the home literacy environment and early number skills was observed (Huntsinger et al., 2016; LeFevre et al., 2009; Segers et al., 2015). The home literacy environment indices used in studies reporting null findings include predominantly meaning-focused literacy items with limited or no code-focused items.

**Why Do Code-focused Home Literacy Experiences Predict Early Number Skills?**

Although numerous studies have established links between home learning experiences and later academic attainment the question of whether such home experiences directly support children's learning remains controversial. Essentially, there are three possible explanations of such relationships. First, more frequent home learning experiences may directly support the development of children's academic skills by providing them with more opportunities to acquire relevant knowledge and to practice skills. Second, home learning experiences may indirectly promote the development of academic skills by influencing the development of language and cognitive abilities that in turn support the development of academic skills. Finally, the relationship may be non-causal with the frequency of home learning experiences having no impact on children's academic skills. Parents of children with more advanced abilities (that may ultimately be genetically influenced) may provide more frequent or more advanced home learning experiences, but the home experiences themselves have no influence on the children's development (see Puglisi et al., 2017 for a discussion).

There were indirect paths from letter-sound interaction experiences via language abilities to all three number skills in Models 2 and 3. Therefore, we can accept the hypothesis that the relationship between letter-sound interaction experiences and all three early number skills is partially mediated by language. It is plausible that letter-sound interaction experiences exert an indirect influence on all three early number skills via the promotion of language. However, it should also be noted that nonverbal abilities accounted for some of the variance explained in the three early number skills by letter-sound interaction experiences. There is a significant (although weaker) path between letter-sound interaction experiences and nonverbal abilities. Given the verbal nature of the interactions and the nonverbal nature of the cognitive skills it would appear unlikely that such path is underpinned by letter-sound interaction experiences supporting the development of nonverbal reasoning and IC. It seems more likely that parents of children with higher nonverbal abilities provide more frequent letter-sound interaction experiences without these interactions exerting a causal influence on nonverbal abilities. This does not preclude a causal relationship between the frequency of letter-sound interaction experiences and language abilities. The path from letter-sound interaction experiences to language was stronger than that to nonverbal abilities. However, the evidence from this study alone cannot discount the relationship between letter-sound interaction experiences and language being correlational rather than causal. The findings from other studies that have examined the relationship between code-focused home literacy experiences and language abilities (particularly phonological awareness) are mixed. Positive relationships (Foy & Mann, 2003; Hamilton et al., 2016; Sénéchal & LeFevre, 2002) and null findings (Hood et al., 2008; Napoli & Purpura, 2018) have been reported. The extent that the specific home experiences captured by the letter-sound interaction experiences scale can promote growth in phonological skills is an area that merits further investigation. To determine unambiguously whether home letter-sound interaction experiences support early number skills via the promotion of language, fully cross-lagged longitudinal studies where letter-sound interaction experiences, language abilities and early number skills are all measured at multiple time points are required.

In addition to the indirect relationships in Model 2, letter-sound interaction experiences directly explained variance in both counting and number transcoding that was not accounted for by either language or nonverbal abilities. These direct relationships cannot simply reflect a parental response to children's language or nonverbal abilities as these skills are already accounted for in the model. Furthermore, letter-sound interaction experiences predicted *growth* in counting in Model 3 and the regression model when the autoregressor is controlled. It is somewhat surprising that discussing letters and sounds could directly support counting. We suggest that discussing the sounds that letters *represent* may support young children's understanding of representational systems more widely. Understanding that letters can represent sounds could support children's appreciation that number words represent quantities. However, if letter-sound interaction experiences exert an influence via strengthening representational understanding it seems inconsistent that they do not predict significant growth in number transcoding in our study (it just missed significance in Model 3). This may in part be due to the relative stability of number transcoding in the time period we studied (Model 3 shows the beta value for number transcoding at $T_1$ predicting number transcoding at $T_3$ to be twice the weight of the equivalent pathway for counting). The extent that letter-sound interaction experiences can predict number transcoding skills merits further investigation.

**Limitations and Further Work**

Although our findings are consistent with letter-sound interaction experiences supporting early number skills, particularly counting, for causal relationships to be unambiguously identified randomised controlled trials of interventions that promote letter-sound interaction experiences in the home are required. For example, such studies could determine whether low intensity interventions promoting parental talk about letter and sounds influences the development of both language and number skills in pre-schoolers.

Possible interactions among early number skills and the home and preschool environments also need to be considered. For example, if preschools provide high-quality literacy and numeracy experiences individual differences in home experiences may have limited impact. Equally, preschool experiences of lower quality may magnify the influence of the child's experiences at home. Multi-centre studies with both large numbers of preschool settings and large numbers of children at each setting to enable multi-level modelling and the systematic exploration of these factors are required. Studies such as Anders et al. (2013) and Melhuish et al. (2008a) have begun to explore these issues, but they need to be extended.

This study extends our understanding of the relationships between the HLE and specific early number skills. Although we acknowledge we have examined only some aspects of numeracy development, counting, number transcoding and calculation are important skills and explained significant independent variance in overall mathematical attainment. To improve further our understanding of the role of the HLE on numeracy development, future studies should examine both non-symbolic number skills (e.g., non-symbolic set comparison) and a wider range of symbolic number skills (e.g., symbolic number comparison, number line estimation). Similarly, our home number experiences scale focused on direct experiences that had an explicit numerical component. This measure could be widened to examine the impact of informal home number experiences (Skwarchuk et al., 2014) or broader mathematical experiences (Hart et al., 2016; Huntsinger et al., 2016; Zippert & Rittle-Johnson, 2020). Although such experiences have not typically been associated with symbolic early number skills (Skwarchuk et al., 2014), they may be associated with non-symbolic number skills or wider mathematical skills. We also acknowledge that although we considered a greater range of language and cognitive abilities than in most HLE studies, our battery of cognitive assessments could be extended further. In particular, we focused on IC whereas future studies could consider executive functioning more comprehensively and also assess shifting and updating of working memory.

**Conclusion**

Code-focused home literacy experiences have been longitudinally associated with emergent literacy skills even when child language and parental skills are controlled (Hood et al., 2008; Puglisi et al., 2017; Sénéchal & LeFevre, 2002). The present study extends these findings by demonstrating that code-focused home literacy experiences longitudinally predict counting and number transcoding skills when both children’s language and nonverbal abilities are controlled. Furthermore, these home experiences predict counting skills even when initial skill levels are included in the model. This suggests that code-focused home literacy experiences could be beneficial in supporting early number skills, particularly counting. We acknowledge that the variance explained by code-focused home literacy experiences, although significant, is modest and that many other factors may influence the development of early number skills. However, given the evidence that code-focused home literacy experiences are beneficial in developing emergent literacy as well as indications that they may support the development of counting skills, parents should have access to information on how to integrate age-appropriate letter-sound interaction experiences into the everyday activities of their young children.

**Table 1**

*Bivariate Correlations and Descriptive Statistics for the Demographic Variables, the Home Learning Scales and the Early Number Skills Tasks*

| | 1 | 2 | 3 | 4 | 5 | 6 | 7 | 8 | 9 | 10 | 11 | 12 | 13 | 14 | 15 | 16 | 17 | 18 | 19 | 20 | 21 |
|---|---|---|---|---|---|---|---|---|---|---|---|---|---|---|---|---|---|---|---|---|---|
| 1. PD | - | | | | | | | | | | | | | | | | | | | | |
| 2. PQ | **.42** | - | | | | | | | | | | | | | | | | | | | |
| 3. N EXP | -.03 | -.01 | - | | | | | | | | | | | | | | | | | | |
| 4. L MEA | .01 | .08 | **.51** | - | | | | | | | | | | | | | | | | | |
| 5. LS INT | .02 | .01 | **.73** | **.50** | - | | | | | | | | | | | | | | | | |
| 6. L ACT | -.11 | **-.16** | **.72** | **.42** | **.70** | - | | | | | | | | | | | | | | | |
| 7. BE | **.19** | **.21** | **.17** | **.14** | **.18** | .03 | - | | | | | | | | | | | | | | |
| 8. SQC $T_1$ | .09 | .01 | **.18** | -.01 | **.22** | .11 | **.13** | - | | | | | | | | | | | | | |
| 9. GMX $T_1$ | **.13** | .08 | **.25** | .06 | **.31** | **.22** | .07 | **.45** | - | | | | | | | | | | | | |
| 10. CO $T_1$ | .02 | .04 | **.24** | .06 | **.25** | **.19** | .11 | **.37** | **.45** | - | | | | | | | | | | | |
| 11. N RECO $T_1$ | .10 | .10 | **.29** | -.01 | **.32** | **.19** | .07 | **.42** | **.58** | **.41** | - | | | | | | | | | | |
| 12. N READ $T_1$ | **.14** | .09 | **.27** | -.06 | **.29** | **.13** | .08 | **.59** | **.58** | **.44** | **.82** | - | | | | | | | | | |
| 13 ADD $T_1$ | .07 | .06 | **.18** | -.02 | **.25** | .12 | .01 | **.45** | **.46** | **.28** | **.55** | **.60** | - | | | | | | | | |
| 14. SUB $T_1$ | .05 | .08 | **.18** | .06 | **.25** | .12 | .06 | **.37** | **.51** | **.32** | **.46** | **.47** | **.54** | - | | | | | | | |
| 15. SQC $T_3$ | .01 | .02 | **.18** | -.05 | **.24** | **.20** | -.06 | **.41** | **.45** | **.32** | **.44** | **.46** | **.35** | **.33** | - | | | | | | |
| 16. GMX $T_3$ | .04 | .05 | **.24** | .06 | **.30** | **.22** | .03 | **.38** | **.49** | **.39** | **.43** | **.48** | **.44** | **.35** | **.53** | - | | | | | |
| 17. CO $T_3$ | -.08 | .06 | **.23** | .02 | **.23** | **.13** | .03 | **.22** | **.28** | **.24** | **.34** | **.37** | **.33** | **.28** | **.39** | **.39** | - | | | | |
| 18. N RECO $T_3$ | .09 | .09 | **.25** | -.01 | **.31** | **.15** | .07 | **.35** | **.47** | **.36** | **.57** | **.59** | **.38** | **.38** | **.56** | **.64** | **.42** | - | | | |
| 19. N READ $T_3$ | .03 | .10 | **.28** | -.01 | **.31** | **.17** | .04 | **.34** | **.46** | **.34** | **.60** | **.62** | **.40** | **.41** | **.63** | **.61** | **.43** | **.76** | - | | |
| 20. ADD $T_3$ | .11 | .10 | **.18** | -.03 | **.24** | .08 | -.05 | **.33** | **.36** | **.33** | **.44** | **.47** | **.45** | **.40** | **.48** | **.55** | **.41** | **.53** | **.51** | - | |
| 21. SUB $T_3$ | **.16** | **.17** | .11 | -.07 | **.18** | **.14** | -.01 | **.39** | **.44** | **.35** | **.45** | **.45** | **.43** | **.45** | **.49** | **.54** | **.41** | **.57** | **.55** | **.66** | - |
| M | 5.42 | 4.75 | 0.00 | 0.00 | 0.00 | 0.00 | 53.41 | 16.57 | 3.17 | 5.14 | 6.41 | 5.07 | 1.69 | 2.23 | 56.91 | 8.29 | 9.78 | 17.57 | 13.39 | 5.83 | 5.37 |
| SD | 3.32 | 2.00 | 0.89 | 0.88 | 0.91 | 0.89 | 21.63 | 14.23 | 2.47 | 2.72 | 5.32 | 3.99 | 2.25 | 2.23 | 38.43 | 3.53 | 4.08 | 5.96 | 4.93 | 3.40 | 3.41 |
| Max | 10 | 8 | - | - | - | - | 100 | - | 15 | 20 | 20 | 20 | 12 | 12 | - | 15 | 20 | 25 | 20 | 12 | 12 |
| *N* | 271 | 270 | 274 | 274 | 274 | 274 | 268 | 274 | 274 | 274 | 274 | 274 | 274 | 274 | 232 | 232 | 232 | 232 | 232 | 232 | 232 |

*Notes.* PD=Postcode decile; PQ=Parental qualification; N EXP=Number experiences; L MEA=Meaning-focused experiences; LS INT=Letter-sound interactions experiences; L ACT=Letter activities; BE=Book exposure; SQC=Sequential counting; GMX=Give me X; CO=Counting objects; N RECO=Numeral recognition; N READ=Numeral reading; ADD=Additions; SUB=Subtractions. Significant correlations in bold. Values above .12 $p<.05$, values above .17 $p<.01$, values above .22 $p<.001$.

**Table 2**

*Descriptive Statistics and Correlations for the Language and Cognitive Measures*

| | 1 | 2 | 3 | 4 | 5 | 6 | 7 | 8 |
|---|---|---|---|---|---|---|---|---|
| 1. Rhyme awareness | **-** | **.48** | **.35** | **.32** | **.15** | **.20** | **.17** | **.26** |
| 2. Alliteration awareness | | **-** | **.33** | **.38** | **.29** | **.19** | .05 | **.25** |
| 3. Naming vocabulary | | | - | **.37** | **.34** | **.33** | **.24** | **.30** |
| 4. Receptive vocabulary | | | | **-** | **.25** | **.32** | **.26** | **.32** |
| 5. Matrices | | | | | **-** | **.26** | **.18** | **.31** |
| 6. Picture similarities | | | | | | **-** | **.26** | **.31** |
| 7. Big/Little Stroop | | | | | | | **-** | **.30** |
| 8. Fish/Shark d’ | | | | | | | | **-** |
| *M* | 4.36 | 3.73 | 127.12 | 16.99 | 57.92 | 92.60 | 75.70 | 1.74 |
| *SD* | 2.71 | 2.80 | 15.09 | 4.77 | 18.39 | 12.36 | 26.71 | 1.12 |
| *N* | 254 | 263 | 265 | 257 | 265 | 256 | 251 | 242[a] |
| α | .83[b] | .84[b] | .73[b] | .89[b] | .75[b] | .79[b] | .77[c] | - |

*Notes.* Significant correlations in bold. Values above .21 $p$<.05, values above .17 $p$<.01, values above .13 $p$<.001. [a] d’ indexes could not be calculated for 8 children for the Fish/Shark measure due to their random response pattern. All other missing data on this and the other measures is due to child absence during the testing schedule. [b] Inter-tem reliability of the standardisation sample. [c] Inter-tem reliability of older preschool sample (Kochanska et al., 2000).

**Table 3**
*Correlations Among the Demographic Variables, the Home Learning Scales and the Number, Cognitive and Language Composite Measures*

| | 1 | 2 | 3 | 4 | 5 | 6 | 7 | 8 | 9 | 10 | 11 | 12 | 13 | 14 | 15 |
|---|---|---|---|---|---|---|---|---|---|---|---|---|---|---|---|
| 1. PD | **-** | | | | | | | | | | | | | | |
| 2. PQ | **.42** | - | | | | | | | | | | | | | |
| 3. N EXP | -.03 | -.01 | - | | | | | | | | | | | | |
| 4. L MEA | .01 | .08 | **.51** | - | | | | | | | | | | | |
| 5. LS INT | .02 | .01 | **.73** | **.50** | - | | | | | | | | | | |
| 6. L ACT | -.11 | **-.16** | **.72** | **.42** | **.70** | - | | | | | | | | | |
| 7. BE | **.19** | **.21** | **.17** | **.14** | **.18** | .03 | - | | | | | | | | |
| 8. COUNT $T_1$ | .09 | .07 | **.29** | .07 | **.33** | **.24** | .11 | - | | | | | | | |
| 9. N TRANSC $T_1$ | **.13** | .10 | **.29** | -.04 | **.32** | **.17** | .08 | **.62** | **-** | | | | | | |
| 10. CALC $T_1$ | .07 | .08 | **.20** | .02 | **.29** | **.13** | .04 | **.53** | **.62** | **-** | | | | | |
| 11. LANG $T_2$ | .06 | .09 | **.20** | .08 | **.31** | **.16** | .12 | **.44** | **.45** | **.44** | **-** | | | | |
| 12. N-VERB $T_2$ | -.06 | -.02 | **.16** | .03 | **.20** | **.21** | .11 | **.39** | **.36** | **.34** | **.49** | **-** | | | |
| 13. COUNT $T_3$ | -.02 | .06 | **.28** | .05 | **.32** | **.25** | .04 | **.49** | **.51** | **.48** | **.45** | **.35** | **-** | | |
| 14. N TRANSC $T_3$ | .06 | .10 | **.28** | -.10 | **.33** | .15 | .06 | **.51** | **.66** | **.48** | **.44** | **.38** | **.67** | **-** | |
| 15. CALC $T_3$ | **.15** | **.15** | **.16** | -.05 | **.23** | **.14** | -.03 | **.47** | **.52** | **.55** | **.46** | **.33** | **.63** | **.63** | **-** |
| *N* | 271 | 270 | 274 | 274 | 274 | 274 | 268 | 274 | 274 | 274 | 253 | 236 | 232 | 232 | 232 |

*Notes.* PD=Postcode decile; PQ=Parental qualification; N EXP=Number experiences; L MEA=Meaning-focused experiences; LS INT=Letter-sound interactions experiences; L ACT=Letter activities; BE=Book exposure; COUNT=Counting; N. TRANSC=Number transcoding; CALC=Calculations; LANG=Language abilities; N-VERB=Nonverbal abilities. Significant correlations in bold. Values above .13 $p$<.05, values above .17 $p$<.01, values above .22 $p$<.001.

**Table 4**

*Hierarchical Regressions Analysing the Relations between Letter-sound Interaction Experiences at $T_1$ and the Number Skills at $T_3$, Controlling for the Autoregressive Effect of Initial Number Skill Level, Language and Nonverbal Abilities*

| | | Counting $T_3$ | Number transcoding $T_3$ | Calculation $T_3$ |
|---|---|---|---|---|
| Model 1 | Step 1. Autoregressor | | | |
| | Preschool number skill ($T_1$) | .49*** | .64*** | .53*** |
| | *F* | 63.59 | 141.84 | 77.41 |
| | $R^2$ | .24*** | .42*** | .28*** |
| Model 2 | Step 1. Autoregressor | | | |
| | Preschool number skill ($T_1$) | .34*** | .54*** | .40*** |
| | Step 2. Language and nonverbal abilities | | | |
| | Language ($T_2$) | .25** | .15* | .23** |
| | Nonverbal abilities ($T_2$) | .11 | .11 | .08 |
| | *F* | 30.75 | 55.11 | 33.61 |
| | $\mathrm{R}^2$ | .32*** | .46*** | .34*** |
| | $\Delta R^2$ | .08*** | .04** | .06*** |
| Model 3 | Step 1. Control Variables | | | |
| | Preschool number skill ($T_1$) | .31*** | .52*** | .40*** |
| | Step 2. Language and nonverbal abilities | | | |
| | Language ($T_2$) | .21** | .13* | .23** |
| | Nonverbal abilities ($T_2$) | .10 | .11 | .08 |
| | Step 3. Home learning | | | |
| | Letter-sound interactions experiences ($T_1$) | .14* | .09 | .00 |
| | *F* | 24.78 | 42.16 | 25.08 |
| | $R^2$ | .34*** | .46*** | .34*** |
| | $\Delta R^2$ | .02* | .01 | .00 |

*Notes*. F values relate to $R^2$. If not specified values reported are the standardised regression coefficients Model 1 *df*=1, 200; Model 2 *df*=3,198; Model 3 *df*=4,197. *$p$<.05, two-tailed. **$p$<.01, two-tailed. ***$p$<.001 two-tailed.

**Figure 1**

*Longitudinal Path Analysis Predicting the Number Skills at Time 3 from the Home Learning Scales and the Demographic Variables at time 1*

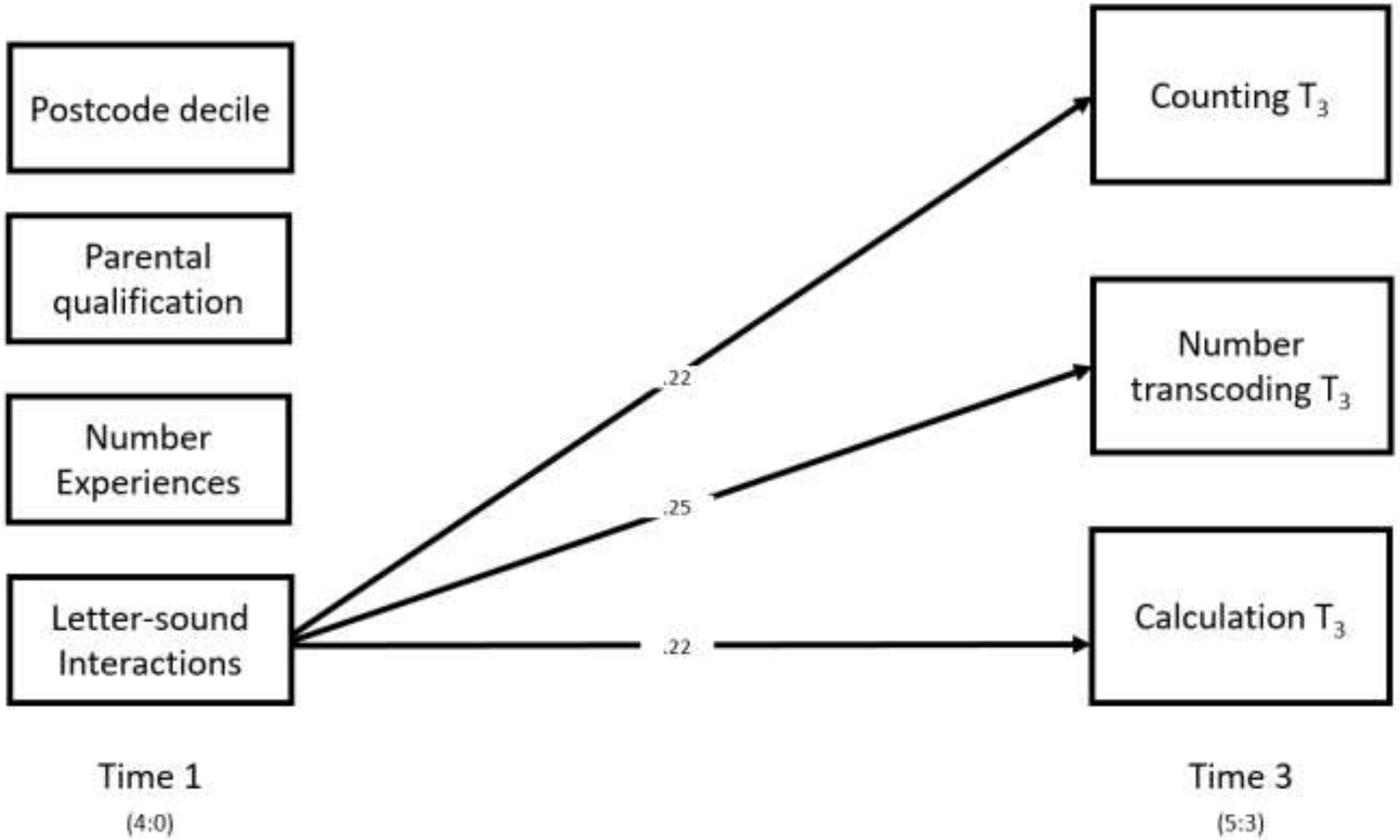

*Note*. Letter-sound interactions = Letter-sound interaction experiences. Only statistically significant paths ($p<.05$) illustrated. $N=226$, cases with missing data are excluded from the presented model. Repeating the model using FIML approach does not result in significant differences. All significant paths remain significant and no other paths become significant.

**Figure 2**

*Longitudinal Path Analysis Model Predicting the Number Skills at Time 3 from the Letter-sound Interactions Scale at Time 1 and the Language and Cognitive Composite Measures at Time 2*

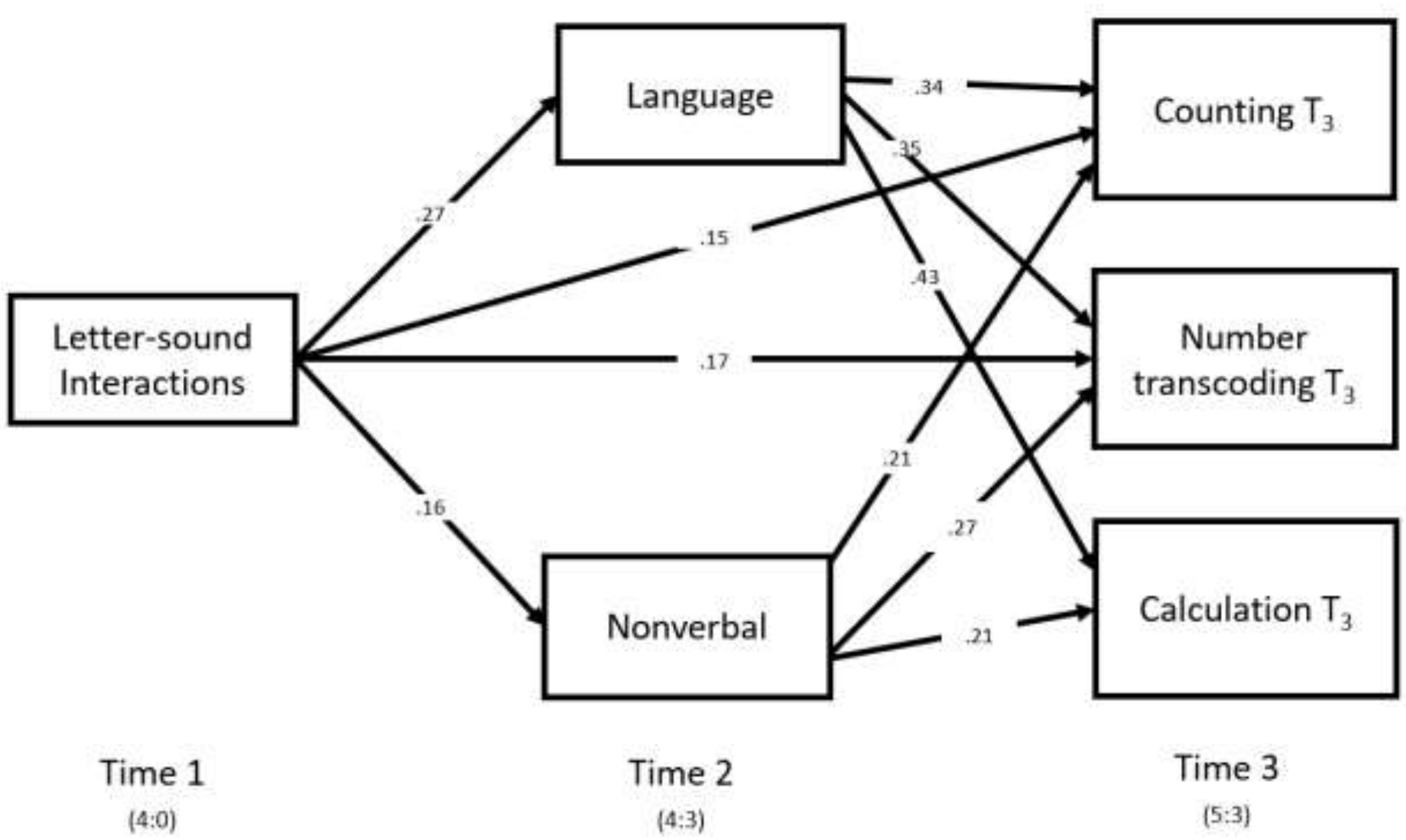

*Note*. Letter-sound interactions = Letter-sound interaction experiences. Only statistically significant paths ($p<.05$) illustrated. $N=204$, cases with missing data are excluded from the presented model. Repeating the model using FIML approach does not result in significant differences. All significant paths remain significant and no other paths become significant.

**Figure 3**

*Longitudinal Path Analysis Model Predicting the Number Skills at Time 3 from the Letter-sound Interactions Scale at Time 1 and the Language and Cognitive Composite Measures at Time 2 with the Autoregressive Effect of Time 1 Number Skills Controlled*

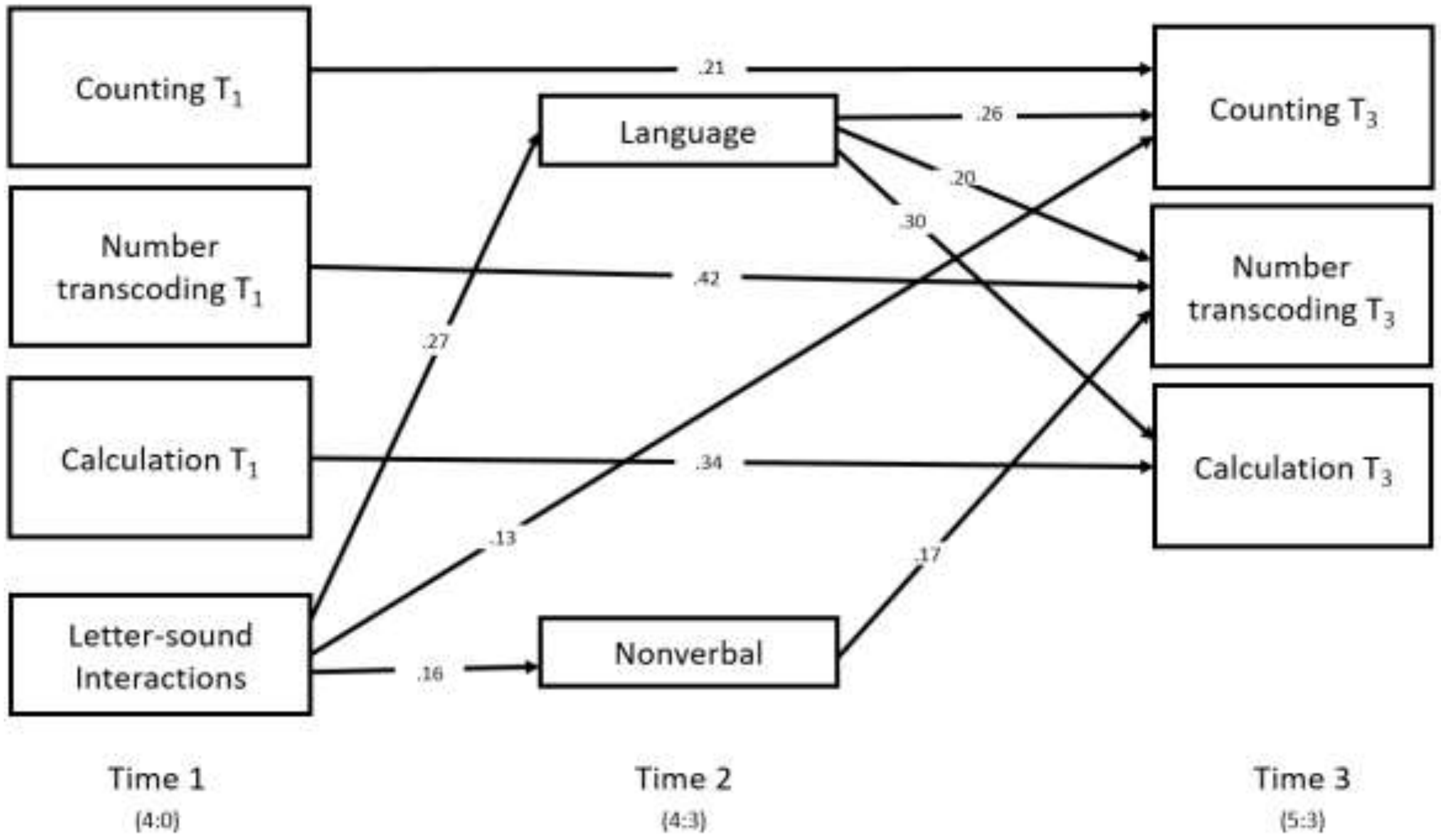

*Note*. Letter-sound interactions = Letter-sound interaction experiences. Only statistically significant paths ($p<.05$) illustrated. Note the poor fit of this model discussed in the results section. *N*=204, cases with missing data are excluded from the presented model. Repeating the model using FIML approach results in minimal differences. All significant paths remain significant. Only the path from nonverbal abilities to counting becomes significant.

**Appendix A**

*Descriptive Statistics and Factor Loadings for Exploratory Factor Analysis Using the Principal Axis Method of the Home Number Experiences Items*

| Item | Mean | SD | Missing Responses | Factor loading |
|---|---|---|---|---|
| Is encouraged to point out or identify numbers in books or the environment (e.g. "What number is on the bus? Can you see a number 8?") | 3.72 | 1.32 | 2 | .74 |
| Is taught the names of numbers (e.g. "This is number 8") | 3.43 | 1.21 | 2 | .67 |
| Writes or traces numbers | 2.32 | 1.46 | 2 | .64 |
| Completes number activities in magazines or workbooks | 1.73 | 1.33 | 3 | .59 |
| Plays games that involve number cards, dice or a number spinner | 2.07 | 1.25 | 1 | .54 |
| Discusses numbers or quantity with an adult (e.g. "How many blocks are there?", "Who has more sandwiches?") | 3.81 | 1.18 | 4 | .51 |
| Recites numbers in order [a] | 4.17 | 1.01 | 1 | |
| Sings number songs (e.g. *Ten Little Monkeys*, *This Old Man*) [b] | 3.47 | 1.26 | 2 | |
| Number experiences [c] | 3.08 | 0.83 | | |

*Notes. N=274,* α=.79. [a] Item excluded from the analyses due to lack of variability in response (operationalized as a mean item score within 1 point of the maximum or minimum or a standard deviation < 1). [b] Item excluded from EFA due to low inter-item correlations within the scale. [c] The mean response to the items contained in the number scale.

**Appendix B**

*Descriptive Statistics and Factor Loadings for Exploratory Factor Analysis Using the Principal Axis Method of the Home Meaning-focused Literacy Experiences Items*

| Item | Mean | SD | Missing Responses | Factor loading |
|---|---|---|---|---|
| Discusses stories with an adult (e.g. "What do you think happens next? Do you think the bunny is frightened?") | 3.70 | 1.08 | 3 | .70 |
| Is encouraged to point out or identify pictures in books (e.g. "Can you point to the elephant?") | 3.88 | 1.07 | 2 | .65 |
| Is encouraged to choose books that interest them to look at with an adult | 3.78 | 1.03 | 1 | .63 |
| Is encouraged to use books to follow-up interests or experiences they have (e.g. looking at a space book because that had talked about space at preschool) | 2.24 | 1.35 | 3 | .58 |
| Discusses with an adult how things work or what they mean (e.g. "Why do you think the ice lolly is melting?", "Nocturnal animals sleep in the day") | 3.62 | 1.42 | 4 | .51 |
| Looks at factual books (e.g. books about animals, space or transport) | 2.90 | 1.27 | 1 | .49 |
| Has stories read to them [a] | 4.14 | 0.78 | 0 | |
| Makes up songs, stories or rhymes [b] | 3.72 | 1.41 | 3 | |
| Meaning-focused literacy experiences [c] | 3.50 | 0.74 | | |

*Notes. N*=274, α=.76. [a] Item excluded from the analyses due to lack of variability in scores (operationalized as a mean item score within 1 point of the maximum or minimum or a standard deviation < 1). [b] Item excluded from EFA due to low inter-item correlations within the scale (operationalized as a correlation coefficient < .3). [c] The mean response to the items contained in the meaning-focused literacy scale.

**Appendix C**

*Descriptive Statistics and Factor Loadings for Exploratory Factor Analysis Using the Principal Axis Method of the Home Code-focused Literacy Experiences Items*

| Item | Mean | SD | Missing Responses | Factor 1 loading | Factor 2 loading |
|---|---|---|---|---|---|
| Is prompted to identify letters in books or the environment (e.g. "Can you see a's' on the sign?", "What letter does the word cat begin with?") | 3.21 | 1.50 | 1 | .82 | |
| Talks about letter sounds with an adult (e.g. "What sound does snake start with?", "Can you think of any other words starting with 's'"? | 3.28 | 1.40 | 2 | .82 | |
| Is taught the names or sounds of letters or how to 'sound out' words | 3.64 | 1.32 | 0 | .52 | .17 |
| Forms or traces letters or writes their name | 2.78 | 1.53 | 3 | .40 | .16 |
| Plays with puzzles or games involving letters | 2.68 | 1.30 | 1 | | .79 |
| Sings or recites the alphabet | 2.93 | 1.51 | 4 | | .67 |
| Completes activities involving letters or sounds in magazines or workbooks | 1.89 | 1.33 | 4 | | .62 |
| Letter-sound interaction experiences [a] | 3.24 | 1.10 | | | |
| Letter activities [b] | 2.50 | 1.14 | | | |

*Notes. N= 274.* [a] The mean response to the items contained in the letter-sound interactions scale (factor 1, α=.76). [b] The mean response to the items contained in the letter activities scale (factor 2, α=.74).

**Appendix D**

*Parents' Responses to the Items in the Book Exposure Measure*

| Item | Identified as Real (%) |
|---|---|
| *Real titles* | |
| The very hungry caterpillar [a] | 251 (91.6) |
| Kipper [b] | 227 (82.8) |
| Dear zoo [a] | 211 (77) |
| That's not my monkey [b] | 202 (73.7) |
| Aliens love underpants [b] | 195 (71.2) |
| The snail and the whale [a] | 193 (70.4) |
| Giraffes can't dance [a] | 172 (62.8) |
| Maisy's bedtime [b] | 161 (58.8) |
| Not now, Bernard [c] | 149 (54.4) |
| Each peach, pear, plum [a] | 129 (47.1) |
| Princess Smartypants [c] | 105 (38.3) |
| Dogger [c] | 89 (32.5) |
| Gorilla [c] | 68 (24.8) |
| Would you rather... [c] | 66 (24.1) |
| Oscar got the blame [b] | 57 (20.8) |
| *Made-up titles* | |
| Grandmother Windmill | 7 (2.6) |
| Belinda Brown takes charge | 15 (5.5) |
| The peg dolly | 19 (6.9) |
| Sally-Anne drives the van | 19 (6.9) |
| What's after bedtime? | 25 (9.1) |
| The wand that wouldn't work | 43 (15.7) |

*Notes.* $n$=268 (six parents did not complete this section of the questionnaire), α=.77 (real titles only). Percentages are provided in brackets. Sources: [a] Amazon best-selling book titles; [b] Most borrowed authors; [c] Booktrust 100 best books.